\documentclass[english, final,leqno,onefignum,onetabnum]{siamltex1213}
\usepackage[T1]{fontenc}
\usepackage[latin9]{inputenc}
\usepackage{color}
\usepackage{amsmath}
\usepackage{amssymb}
\usepackage{setspace}
\usepackage{graphics, epsfig}

\usepackage{algorithm}
\usepackage{algorithmicx}
\usepackage{algpseudocode}

\begin{document}
\title{On the Convergence Rate of Variants of the Conjugate Gradient Algorithm
in Finite Precision Arithmetic}

\author{Anne Greenbaum%
\and
Hexuan Liu%
\and
Tyler Chen%
    \thanks{University of Washington, Applied Mathematics Dept., Box 353925,
            Seattle, WA 98195.  This work was supported in part by NSF grant DMS-1210886.}
}
\maketitle

\begin{abstract}
We consider three mathematically equivalent variants of the conjugate gradient
(CG) algorithm and how they perform in finite precision arithmetic.  
It was
shown in [{\em Behavior of slightly perturbed Lanczos and conjugate-gradient
recurrences}, Lin.~Alg.~Appl., 113 (1989), pp.~7-63] that under certain
conditions 
the convergence of a slightly perturbed CG computation
is like that
of exact CG for a matrix with many eigenvalues distributed throughout tiny
intervals about the eigenvalues of the given matrix, the size
of the intervals being determined by how closely these conditions are satisfied. 
We determine to what
extent each of these variants satisfies the desired conditions, using a set
of test problems and show that there is significant correlation
between how well these conditions are satisfied and how well the finite
precision computation converges before reaching its ultimately attainable
accuracy.  
We show that for problems where the width of the intervals containing 
the eigenvalues
of the associated exact CG matrix makes a
significant difference in the behavior of exact CG, the different CG variants
behave differently in finite precision arithmetic.
For problems where the interval width makes little difference or where the
convergence of exact CG is essentially governed by the upper bound based
on the square root of the condition number of the matrix, the different CG
variants converge similarly in finite precision arithmetic until the ultimate
level of accuracy is achieved, although this ultimate level of accuracy may
be different for the different variants.  This points to the
need for testing new CG variants on problems that are especially sensitive
to rounding errors.
\end{abstract}

\section{Introduction}
Several variants of the conjugate gradient algorithm (CG) for solving a symmetric positive
definite linear system $A x = b$ have been proposed to make better use of parallelism; 
see, e.g.,  \cite{Ros83,Saa85,Meu87,Saa89,ChronGear,ChronGear2,GhyVan}.  
Here we consider three variants:
the original Hestenes and Stiefel algorithm \cite{HS} (HSCG), a variant with
somewhat more opportunity for parallelism due to Chronopoulos and Gear \cite{ChronGear} (CGCG),
and a more recent pipelined version due to Ghysels and Vanroose \cite{GhyVan} (GVCG).
While all of these algorithms are mathematically equivalent, they behave differently 
when implemented in finite precision arithmetic.  Perhaps the most
dramatic difference is in the ultimately 
attainable accuracy of the computed solution.  All of these algorithms compute
an initial residual $r_0 = b - A x_0$, where $x_0$ is the initial guess for the solution,
and then compute updated ``residual'' vectors $r_k$, $k=1,2, \ldots$, using a recurrence
formula.  In finite precision arithmetic, however, these updated vectors may differ
from the true residuals $b - A x_k$, where $x_k$ is the approximate solution vector
generated at step $k$.  When this difference becomes large, 
the norms of the updated vectors may or may not continue to decrease, but the true 
residual norm (that is, the norm of $b - A x_k$) levels off (or may 
even grow).
The level of accuracy of the approximate solution $x_k$ when this occurs is studied 
in \cite{Carson,Cools}.

In this paper, we consider what happens in the stage
{\em before} the true and updated
residual vectors start to deviate significantly.  Even during this stage, the 
different variants may show different convergence patterns on problems with
certain eigenvalue distributions that make them especially sensitive to rounding errors.
This is a phenomenon that we wish to understand.  On other problems, 
where eigenvalues of the coefficient matrix are distributed in a more uniform way,
the algorithms may all behave very similarly.  This, too, is something that needs
a mathematical explanation since this may hold {\em even after agreement with exact arithmetic
is lost}.  Specifically, we aim to give an explanation for the
results shown in Figure \ref{fig:bcsstk03}, where different CG variants 
converge differently, and in Figure \ref{fig:bcsstkall}, where, although one
variant fails to obtain as accurate a solution as the others, up until the point where the
convergence curve of that computation levels off, all variants converge similarly.

A good deal of work beginning in the 1980's (and in the thesis of Paige \cite{Paige2}
dating back to 1971) has been aimed at explaining the behavior
of the Lanczos and conjugate gradient algorithms in finite precision arithmetic,
or, more generally, when the recurrence formulas for these algorithms
are perturbed slightly.  See, for example, \cite{Paige1,Greenbaum89,GreenStrak,DGK}.  
In exact arithmetic, the Lanczos algorithm can be thought of as a
part of the CG algorithm:  The Lanczos algorithm generates a sequence of tridiagonal matrices and the CG 
algorithm implicitly solves linear systems with these tridiagonal matrices in order
to approximate the solution of the linear system $Ax=b$.
In a seminal paper \cite{Paige1},
Paige showed that what is now the standard implementation of the Lanczos algorithm 
maintained certain local orthogonality and normalization properties
even when implemented in finite precision arithmetic and that those properties
could be used to establish results about the eigenvalue/vector approximations
generated during a finite precision Lanczos computation.  A nice summary of this work
can be found in \cite{Parlett}, and more recent work by Paige has
extended these results significantly \cite{Paige2010,Paigetoappear}.  
Later, these same properties were
used in \cite{Greenbaum89} to establish results about the 
convergence of the conjugate gradient 
algorithm under the {\em assumption} that these properties were satisfied by a finite
precision, or otherwise slightly perturbed, implementation.  
A natural question to ask is:  Which of the various proposed
implementations satisfy these properties, and do those that do satisfy the properties
used in Paige's analysis have better behavior than those that do not?  For those that do not, 
are there other ways to explain their behavior?  
In this paper, we consider three mathematically equivalent variants of the conjugate
gradient algorithm:  the original Hestenes and Stiefel variant \cite{HS} (HSCG), a variant with
somewhat more opportunity for parallelism due to Chronopoulos and Gear \cite{ChronGear} (CGCG),
and a more recent pipelined variant due to Ghysels and Vanroose \cite{GhyVan} (GVCG). 

In the following subsections, we review the properties
that have been assumed in order to relate finite precision or
otherwise slightly perturbed CG computations to exact CG for a larger matrix
with eigenvalues in small intervals about those of the given matrix.
We do not rigorously prove these properties for
any particular implementation but give an indication why one of them may be expected to fail
in the GVCG variant, and we check numerically whether or not they hold for a 
number of test problems and whether satisfaction of such properties coincides with faster
convergence (in terms of number of iterations).  Again, the analysis of this paper
deals only with steps at which the true and updated residuals are still in
close agreement, that is, {\em before} the best level of accuracy is reached.
The reader is referred to \cite{Carson,Cools} for a discussion of what this best level 
of accuracy is for different CG variants.

Throughout the paper, $A$ will denote a real symmetric positive definite matrix, although
the results are easily extended to complex Hermitian positive definite matrices.  The symbol
$\| \cdot \|$ will denote the 2-norm for vectors and the corresponding spectral norm for 
matrices. 

\subsection{Slightly Perturbed Lanczos Computations}
In \cite{Greenbaum89}, an analogy was established between finite precision, 
or otherwise slightly perturbed, Lanczos
computations with matrix $A$ and initial vector $q_1$ and exact Lanczos
applied to a larger matrix $T$ whose eigenvalues all lie in tiny intervals
about the eigenvalues of $A$.  
More precisely, it was shown that if $T_J$ is the tridiagonal matrix
produced at step $J$ of a finite precision computation with matrix $A$, 
then $T_J$ can be extended
to a larger tridiagonal matrix $T$ whose eigenvalues are all close to eigenvalues of $A$,
assuming that the finite precision computation satisfies certain local orthogonality
and normalization properties.  
An algorithm was given for producing such an extension, 
and this algorithm is outlined in the Appendix of this paper.
When exact Lanczos is applied to $T$ with initial vector $\hat{q}_1$ equal to the first unit vector,
it produces the same tridiagonal matrices $T_1 , \ldots T_J$ as the finite precision 
computation.  If $\lambda_{i \ell}$, $\ell =1,2, \ldots$, are the eigenvalues of $T$ that
are close to eigenvalue $\lambda_i$ of $A$, and if $v_{i \ell}$, $\ell = 1,2, \ldots$,
are the corresponding orthonormal eigenvectors of $T$, while $v_i$ is a normalized 
eigenvector of $A$ associated with $\lambda_i$, then it was shown that
\[
\sum_{\ell} \langle \hat{q}_1 , v_{i \ell} \rangle^2 \approx \langle q_1 , v_i \rangle^2 .
\]
This meant that theorems (that
assume exact arithmetic) about the behavior of the first $J$ steps of the
Lanczos algorithm applied to 
such matrices $T$ with such initial vectors $\hat{q}_1$ could be applied to finite
precision computations with matrix $A$ and initial vector $q_1$.  

The assumptions needed for this analogy to hold were that vectors $q_1 , \ldots , q_J$ 
generated by the finite precision computation satisfied
\begin{equation}
A Q_J = Q_J T_J + \beta_J q_{J+1} \xi_J^T + F_J , \label{FPLan}
\end{equation}
where $Q_J$ is the $n$ by $J$ matrix whose columns are $q_1 , \ldots , q_J$, 
$T_J$ is the symmetric tridiagonal matrix
\[
T_J = \left[ \begin{array}{ccccc}
\alpha_1 & \beta_1 & & & \\
\beta_1  & \ddots  & \ddots & & \\
         & \ddots  & \ddots & \beta_{J-1} & \\
         &         & \beta_{J-1} & \alpha_J \end{array} \right] ,
\]
$\xi_J$ is the $J$th unit vector $(0, \ldots , 0, 1 )^T$, and $F_J$, which accounts for
rounding errors, has columns $f_j$, $j=1, \ldots , J$, satisfying
\begin{equation}
\| f_j \| \leq \epsilon \| A \| , \label{FJnorm}
\end{equation}
where $\epsilon$ is tiny (ideally, a modest multiple of the machine precision).  
It is further assumed that,
because of the choice of the coefficients, $\alpha_1 , \ldots , \alpha_J , \beta_1 , \ldots , 
\beta_J$, the 2-norm of each vector $q_j$ is approximately $1$, say, 
$\| q_j \| \in [1- \epsilon , 1+ \epsilon]$,
and the inner product of {\em successive} pairs of vectors is almost $0$:
\begin{equation}
| \langle \beta_j q_{j+1} , q_j \rangle | \leq \epsilon \| A \| . \label{successiveinnerprod}  
\end{equation}
The analysis in \cite{Greenbaum89}
applies to any computation that satisfies these assumptions for some $\epsilon \ll 1$.
One can extend the tridiagonal matrix $T_J$ in the way described in \cite{Greenbaum89}
even if $\epsilon$ is not very small, but then there is no guarantee that the eigenvalues of
the extended matrix $T$ will all be close to those of $A$ or even that $T$ will be positive definite.
This analysis relies heavily on the work of Paige \cite{Paige1,Paige2},
who showed that a good finite precision implementation of the Lanczos algorithm
satisfies these assumptions, with explicit bounds on the quantities denoted here
as $\epsilon$.

\subsection{Relation Between CG Residuals and Lanczos Vectors}
One way to solve a symmetric positive definite linear system $Ax=b$ is
to use the Lanczos algorithm to generate the matrices $Q_J$ and $T_J$ in (\ref{FPLan}),
taking $q_1$ to be the normalized initial residual:  If $x_0$ is an initial guess for
the solution of the linear system and $r_0 = b - A x_0$, then $q_1 = r_0 / \beta$, 
where $\beta = \| r_0 \|$.  Then the approximate solution $x_J$ is taken to be 
\begin{equation}
x_J = x_0 + Q_J T_J^{-1} \beta \xi_1 . \label{xJLan}
\end{equation}
This choice of $x_J$ minimizes the $A$-norm of the error, $\langle A^{-1} b - x_J , b - A x_J \rangle^{1/2}$,
among all vectors of the form $x_0 + Q_J y$, $y \in \mathbb{R}^J$.  See, for example,
\cite[Sec.~2.5]{Greenbaumbook}.  The conjugate gradient algorithm does this implicitly by generating
Cholesky-like factorizations of the successive tridiagonal matrices $T_k$, $k=1, \ldots , J$:
$T_k = L_k D_k L_k^T$, where $L_k$ is a unit lower bidiagonal matrix and $D_k$ is a diagonal 
matrix.  Thus, rounding errors in the conjugate gradient algorithm involve not only those
in constructing the columns of $Q_J$ but also those in solving the linear system involving
$T_J$.  Since the two processes are intertwined, the effect of rounding errors can be
more difficult to analyze in the conjugate gradient algorithm than in the Lanczos algorithm.
 
The conjugate gradient algorithm for solving a symmetric positive definite linear system
$Ax=b$ can be written in the following form, due to Hestenes and Stiefel \cite{HS} (HSCG):  

\begin{center}
\begin{minipage}{.9\linewidth}
\begin{algorithm}[H]
\caption{Hestenes and Stiefel Conjugate Gradient}\label{hscg_alg}
\fontsize{10}{10}\selectfont
\begin{algorithmic}
\Procedure{HSCG}{$A$, $b$, $x_0$};~~~{\bf Output}($x_k$)
    \State \textbf{Set} $r_0 = b-Ax_0$, $\nu_0 = \langle r_0,r_0 \rangle$, $p_0 = r_0$, $s_0 = Ar_0$, $a_0 = \nu_0 / \langle p_0,s_0 \rangle$.
\For {$k=1,2,\ldots$}
    \State $x_k = x_{k-1} + a_{k-1} p_{k-1}$
    \State $r_k = r_{k-1} - a_{k-1} s_{k-1}$
    \State $\nu_{k} = \langle r_k,r_k \rangle$, $~~b_k = \nu_k / \nu_{k-1}$
    \State $p_k = r_k + b_k p_{k-1}$
    \State $s_k = A p_k$
    \State $\mu_k = \langle p_k,s_k \rangle$, $~~a_k = \nu_k / \mu_k$ 
   \EndFor
\EndProcedure
\end{algorithmic}
\end{algorithm}
\end{minipage}
\end{center}
\vspace{.2in}

As noted above,
it is well-known that if $q_1 = r_0 / \| r_0 \|$ in the Lanczos algorithm, then
subsequent Lanczos vectors are just scaled versions of the corresponding CG
residuals.  To see this from the HSCG algorithm, we first note that the residual vectors 
$r_k$, $k=0,1, \ldots $, satisfy a 3-term recurrence:
\begin{eqnarray*}
r_k & = & r_{k-1} - a_{k-1} A p_{k-1} = r_{k-1} - a_{k-1} A ( r_{k-1} + b_{k-1} p_{k-2} ) \\
    & = & r_{k-1} - a_{k-1} A r_{k-1} - \frac{a_{k-1} b_{k-1}}{a_{k-2}} ( r_{k-2} - r_{k-1} ) . 
\end{eqnarray*}  
If we define normalized residuals by $q_{k+1} := (-1 )^k \frac{r_k}{\| r_k \|}$, then
these vectors satisfy
\[
q_{k+1} = a_{k-1} \frac{\| r_{k-1} \|}{\| r_k \|} A q_k - \left( 1 + \frac{a_{k-1} b_{k-1}}
{a_{k-2}} \right) \frac{\| r_{k-1} \|}{\| r_k \|} q_k - 
\frac{a_{k-1} b_{k-1}}{a_{k-2}} \frac{\| r_{k-2} \|}{\| r_k \|} q_{k-1} ,
\]
or,
\[
A q_k = \frac{\| r_k \|}{ a_{k-1} \| r_{k-1} \|} q_{k+1} + \left( \frac{1}{a_{k-1}} +
\frac{b_{k-1}}{a_{k-2}} \right) q_k + \frac{b_{k-1}}{a_{k-2}} \frac{\| r_{k-2} \|}{\| r_{k-1} \|}
q_{k-1} .
\]
Finally, noting that $b_{k-1} = \| r_{k-1} \|^2 / \| r_{k-2} \|^2$, this becomes
\begin{equation}
A q_k = \frac{\| r_k \|}{ a_{k-1} \| r_{k-1} \|} q_{k+1} + \left( \frac{1}{a_{k-1}} +
\frac{b_{k-1}}{a_{k-2}} \right) q_k + \frac{\| r_{k-1} \|}{a_{k-2} \| r_{k-2} \|} q_{k-1} .
\label{CGLan}
\end{equation}
Thus, if $Q_J$ is the $n$ by $J$ matrix whose columns are $q_1 , \ldots , q_J$, then
\begin{equation}
A Q_J = Q_J T_J + \beta_J q_{J+1} \xi_J^T , \label{CGLanMatrix}
\end{equation}
where $\beta_J = \| r_J \| / ( a_{J-1} \| r_{J-1} \| )$ and $T_J$ is a symmetric tridiagonal
matrix with diagonal entries $\alpha_j = 1/ a_{j-1} + b_{j-1} / a_{j-2}$, $j=1,2, \ldots , J$,
(where terms involving $a_{-1}$ are taken to be $0$) and sub and super diagonal entries 
$\beta_j = \| r_j \| / ( a_{j-1} \| r_{j-1} \| )$, $j=1, \ldots , J-1$.
It follows that if formula (\ref{CGLanMatrix}) can be replaced by something of the form
(\ref{FPLan}) when the columns of $Q_J$ come from normalizing ``residual'' vectors $r_k$ 
in a finite precision CG computation, with the computed vectors satisfying properties
(\ref{FJnorm}) and (\ref{successiveinnerprod}) as well, then
the analysis of \cite{Greenbaum89} will apply to the finite precision CG computation.
We emphasize again that the analysis in \cite{Greenbaum89} gives
information about the rate at which the 
{\em updated} residual vectors $r_k$ decrease in norm and thus is of interest only
as long as these updated vectors resemble the true residuals, $b - A x_k$.  

\subsection{Implications for Finite Precision CG Implementations}
Under the assumption that formulas (\ref{FPLan}), (\ref{FJnorm}), 
(\ref{successiveinnerprod}) hold for some small value $\epsilon$, when the columns of $Q_J$ 
come from normalizing updated CG residual vectors and the entries of $T_J$ are derived from CG coefficients
as described above, the analysis in \cite{Greenbaum89} shows that the updated
CG residual vectors converge at the rate predicted by exact arithmetic theory for a 
symmetric positive definite matrix $T$ whose condition number $\kappa$ is just slightly 
larger than that of $A$:
\begin{equation}
\frac{\| r_k \|}{\| r_0 \|} \leq~ \kappa^{1/2}~ 2 \left( \frac{\sqrt{\kappa} -1}
{\sqrt{\kappa} +1} \right)^k , \label{FPresnorm}
\end{equation} 
although this may be a substantial overestimate.
It shows further that the $A$-norm of the error in the finite precision computation -- that is,
the quantity $\langle r_k , A^{-1} r_k \rangle^{1/2}$ -- is reduced at about the same rate as
the $T$-norm of the error in exact CG applied to $T$:
\begin{equation}
\frac{\langle r_k , A^{-1} r_k \rangle^{1/2}}{\langle r_0 , A^{-1} r_0 \rangle^{1/2}}
\stackrel{<}{\sim} 2 \left( \frac{\sqrt{\kappa} -1}{\sqrt{\kappa} +1} \right)^k , \label{FPAnorm}
\end{equation}
which again may be an overestimate.

A sharper bound on the quantities on the left in (\ref{FPresnorm}) and (\ref{FPAnorm}) 
can be given in
terms of the size of the $k$th degree minimax polynomial on the union of tiny intervals containing the
eigenvalues of $T$; if these intervals are $[ \lambda_i - \delta , \lambda_i + \delta ]$,
then the quantity 
\begin{equation}
\min_{\{ p_k : p_k (0) = 1 \}} \max_{z \in \cup_{i=1}^n 
[ \lambda_i - \delta , \lambda_i + \delta ]} | p_k (z) | \label{sharpbound}
\end{equation}
is an upper bound for
the quantity on the left in (\ref{FPAnorm}) and $\kappa^{1/2}$ times this value is an
upper bound for the left-hand side of (\ref{FPresnorm}).  For some eigenvalue 
distributions, such as eigenvalues fairly uniformly distributed between
$\lambda_{\text{min}}$ and $\lambda_{\text{max}}$, the size of this minimax polynomial is 
not much less than that of the
Chebyshev polynomial on the entire interval $[ \lambda_{\text{min}} - \delta , \lambda_{\text{max}}
+ \delta ]$, on which the bounds in (\ref{FPresnorm}) and (\ref{FPAnorm}) are based.
However, for other
eigenvalue distributions such as eigenvalues tightly clustered at the lower 
end and highly spread out at the upper end, as will be seen in one of our examples, the difference can be great.

These bounds are independent of the initial residual $r_0$.  With knowledge of the size of
components of $r_0$ in the directions of each eigenvector of $A$, the analysis in \cite{Greenbaum89}
gives additional insight into the convergence of a finite precision CG computation that satisfies
the assumptions in \cite{Greenbaum89}.  It behaves
like exact CG applied to a matrix whose eigenvalues lie in tiny intervals about the eigenvalues
of $A$, with an initial residual $\hat{r}_0$ satisfying
\begin{equation}
\sum_{\ell} \langle \hat{r}_0 , v_{i_{\ell}} \rangle^2 \approx \langle r_0 , v_i \rangle^2 ,~~
i=1, \ldots , n , \label{r0hat}
\end{equation}
where $v_i$ is a normalized eigenvector of $A$ corresponding to eigenvalue $\lambda_i$, 
$v_{i_{\ell}}$ is a normalized eigenvector of $T$ corresponding to eigenvalue 
$\lambda_{i_{\ell}}$, and the sum over $\ell$ is the sum over all eigenvalues of $T$ 
in the small interval $[ \lambda_i - \delta , \lambda_i + \delta ]$.  
In some cases, 
even assuming exact arithmetic where
$\delta = 0$, bounds based on the size of the minimax polynomial on the set of eigenvalues
are large overestimates for observed convergence rates. 
While for any given $k$, there
is an initial residual for which equality will hold at step $k$ \cite{Greenbaum79}, components
of that initial residual may differ by hundreds of orders of magnitude.
Such an initial residual
could not even be represented on a machine with standard bounds on exponent size, so whatever the
initial residual in the finite precision computation, it is necessarily far from the
worst possible one. 

\section{CG Variants Designed for Parallelism}

While individual matrix-vector multiplications can be parallelized and vectors can be
partitioned among different processors in the HSCG algorithm, almost none
of the high-level operations comprising an iteration
in that algorithm can be performed simultaneously.
Looking at the algorithm of the previous section, it can be seen that
at each iteration, the matrix-vector product $A p_{k-1}$ must be started,
with at least part of it completed,
before computation of the inner product $\langle p_{k-1} , A p_{k-1} \rangle$ can begin.
This inner product must be {\em completed} before the vectors $x_k$ and $r_k$ can
be formed, and $r_k$ must be at least partly completed
before computation of the next inner product $\langle r_k , r_k \rangle$ can begin.  
This inner product must be {\em completed} before $p_k$ can be formed, and at
least part of $p_k$ must be completed before the start of the next iteration
computing $A p_k$.  It has been observed that waiting for the two
inner products to complete can be very costly when using large numbers of
processors \cite{AGHV12,Cools}.

Several mathematically equivalent CG variants have been devised to allow
overlapping of inner products with each other and with the matrix-vector multiplication 
in each iteration of the 
algorithm.  In the following sections, we consider two of these:  one
due to Chronopoulos and Gear \cite{ChronGear,ChronGear2} (CGCG) that allows either overlapping
of the two inner products or overlapping of one of these with the matrix-vector product,
and a pipelined version due to Ghysels and Vanroose \cite{GhyVan} (GVCG)
that allows overlapping of both inner products as well as the matrix-vector 
multiplication.  We give an indication of why the value of $\epsilon$ in (\ref{FJnorm})
might be expected to be larger for the GVCG variant than for HSCG and CGCG,
and we demonstrate that it is, indeed, larger for a set of test problems.  This does not 
{\em necessarily} mean slower convergence, however; it means simply that in the analogy with exact CG 
for a matrix whose eigenvalues are in intervals about the eigenvalues of $A$, the interval size must be 
larger for GVCG.  If this interval size makes a significant difference in the convergence rate of exact CG,
then we expect the finite precision GVCG computation to require more iterations than
the other variants.  Again, we emphasize that our analysis applies only {\em before}
agreement between true and updated residual vectors is lost and {\em before} the ultimate 
level of accuracy is achieved.  All of our experiments are
performed on a single processor using standard double precision arithmetic,
and we do not consider the timing of the algorithms,
only the number of iterations required to reach a given level of accuracy 
(assuming that that level of accuracy is reached by the algorithm).

Since it is known that a good implementation of the Lanczos algorithm
satisfies (\ref{FPLan} -- \ref{successiveinnerprod}) \cite{Paige1,Paige2}, we also compared 
the CG variants with results obtained by using the Lanczos algorithm to generate 
$Q_J$ and $T_J$ in (\ref{CGLanMatrix}) and then using extra precision to compute
$x_J$ in (\ref{xJLan}).  These results were very similar to those obtained with HSCG,
so they are not included in the plots.  

\subsection{HSCG}\label{subsection:HSCG}
When the Hestenes and Stiefel algorithm of the previous section is implemented in
finite precision arithmetic, the vectors $r_k$ and $p_k$ satisfy
\begin{eqnarray*}
r_k & = & r_{k-1} - a_{k-1} A p_{k-1} + \delta_{r_k} , \\
p_k & = & r_k + b_{k} p_{k-1} + \delta_{p_k} .
\end{eqnarray*}
The roundoff terms can be bounded, as in \cite{StrakosTichy},
using standard results on floating point arithmetic; see, e.g., \cite{Higham}.  For a scalar $\alpha$,
$n$-vectors $v$ and $w$, and an $n$ by $n$ matrix $A$, we have
\begin{eqnarray}
\| \alpha v - \mbox{fl}[ \alpha v ] \| & \leq & \epsilon \| \alpha v \| , \nonumber \\
\| v + w - \mbox{fl}[ v + w ] \| & \leq & \epsilon ( \| v \| + \| w \| ) , \nonumber \\
\| A v - \mbox{fl}[ Av ] \| & \leq & \epsilon c \| A \| \| v \| , \label{standard}
\end{eqnarray}
where $\epsilon$ is the machine precision and $\mbox{fl}[ \cdot ]$ denotes the 
result of floating point evaluation.  If $A$ has at most $m$ nonzeros per row
and the matrix-vector product is computed in the standard way, then $c$
can be taken to be $m n^{1/2}$; alternatively, $c \| A \|$ in (\ref{standard})
can be replaced by $m \|~|A|~\|$, where $|A|$ is the $n$ by $n$ matrix whose
entries are the absolute values of those of $A$.  Applying these rules to
the formulas in the HSCG algorithm, we can write
\begin{eqnarray*}
\| \delta_{r_k} \| & \leq & \epsilon~( \| r_{k-1} \| + 2 \| a_{k-1} A p_{k-1} \| + c
\| A \|~\| a_{k-1} p_{k-1} \| ) + O( \epsilon^2 ) , \\
\| \delta_{p_k} \| & \leq & \epsilon~( \| r_k \| + 2 \| b_k p_{k-1} \| ) 
+ O ( \epsilon^2 ) ,
\end{eqnarray*}
Similarly, the errors in the computed coefficients $a_{k-1}$ and $b_k$ can be
bounded as in \cite{StrakosTichy}.
Here we will assume that the coefficients satisfy the 
formulas in the HSCG algorithm, namely,
\[
a_{k-1} = \frac{\langle r_{k-1} , r_{k-1} \rangle}{\langle p_{k-1} ,A p_{k-1} \rangle} ,~~
b_k = \frac{\langle r_k , r_k \rangle}{\langle r_{k-1} , r_{k-1} \rangle} ,
\]
and we will include any errors in computing these formulas in the $\delta_{r_k}$
and $\delta_{p_k}$ terms.
It follows, as in \cite{Greenbaum89},  that
\begin{eqnarray*}
r_k & = & r_{k-1} - a_{k-1} A ( r_{k-1} + b_{k-1} p_{k-2} + \delta_{p_{k-1}} ) + \delta_{r_k} \\
    & = & r_{k-1} - a_{k-1} A r_{k-1} - \frac{a_{k-1} b_{k-1}}{a_{k-2}} ( r_{k-2} - r_{k-1} 
+ \delta_{r_{k-1}} ) - a_{k-1} A \delta_{p_{k-1}} + \delta_{r_k} \\
    & = & r_{k-1} - a_{k-1} A r_{k-1} - \frac{a_{k-1} b_{k-1}}{a_{k-2}} ( r_{k-2} - r_{k-1} )
- \gamma_k ,
\end{eqnarray*}
where
\begin{equation}
\gamma_k = \frac{a_{k-1} b_{k-1}}{a_{k-2}} \delta_{r_{k-1}} + a_{k-1} A \delta_{p_{k-1}} - \delta_{r_k} .
\label{HSCGgammak}
\end{equation}
Defining $q_{k+1} := (-1 )^k r_k / \| r_k \|$, we can write
\[
q_{k+1} = a_{k-1} \frac{\| r_{k-1} \|}{\| r_k \|} A q_k - \left( 1 + \frac{a_{k-1} b_{k-1}}
{a_{k-2}} \right) \frac{\| r_{k-1} \|}{\| r_k \|} q_k - 
\frac{a_{k-1} b_{k-1}}{a_{k-2}} \frac{\| r_{k-2} \|}{\| r_k \|} q_{k-1} -
\frac{ ( -1 )^k}{\| r_k \|} \gamma_k ,
\]
or,
\[
A q_k = \frac{\| r_k \|}{ a_{k-1} \| r_{k-1} \|} q_{k+1} + \left( \frac{1}{a_{k-1}} +
\frac{b_{k-1}}{a_{k-2}} \right) q_k + \frac{b_{k-1}}{a_{k-2}} \frac{\| r_{k-2} \|}{\| r_{k-1} \|}
q_{k-1} + \frac{(-1 )^k}{a_{k-1} \| r_{k-1} \|} \gamma_k .
\]
With the formula for $b_{k-1}$, this takes a form similar to (\ref{CGLan}):
\begin{equation}
A q_k = \frac{\| r_k \|}{ a_{k-1} \| r_{k-1} \|} q_{k+1} + \left( \frac{1}{a_{k-1}} +
\frac{b_{k-1}}{a_{k-2}} \right) q_k + \frac{\| r_{k-1} \|}{a_{k-2} \| r_{k-2} \|}
q_{k-1} + \frac{(-1 )^k}{a_{k-1} \| r_{k-1} \|} \gamma_k . \label{HSCGLan1}
\end{equation}
From (\ref{HSCGgammak}), the roundoff term in (\ref{HSCGLan1}) can be written as
\begin{equation}
\frac{(-1 )^k}{a_{k-1} \| r_{k-1} \|} \gamma_k =
\frac{( -1 )^k}{\| r_{k-1} \|} \left( \frac{b_{k-1}}{a_{k-2}} \delta_{r_{k-1}} + A \delta_{p_{k-1}} -
\frac{1}{a_{k-1}} \delta_{r_k} \right) . \label{fkHSCG}
\end{equation}
This is the $k$th column of $F_J$ in (\ref{FPLan}).  Note that it 
involves only local rounding errors and might therefore be expected to be of order
$\epsilon \| A \|$, where $\epsilon$ is a moderate multiple of the machine precision.
We will see in later examples that this is indeed the case.


\subsection{CGCG}\label{subsection:CGCG}
Chronopoulos and Gear proposed the following version of the CG algorithm to make better
use of parallelism \cite{ChronGear}:

\begin{center}
\begin{minipage}{.9\linewidth}
\begin{algorithm}[H]
\caption{Chronopoulos and Gear Conjugate Gradient}\label{cgcg_alg}
\fontsize{10}{10}\selectfont
\begin{algorithmic}
\Procedure{CGCG}{$A$, $b$, $x_0$};~~~{\bf Output}($x_k$)
    \State \textbf{Set} $r_0 = b-Ax_0$, $\nu_0 = \langle r_0,r_0 \rangle$, $p_0 = r_0$, $s_0 = Ap_0$, $a_0 = \nu_0 / \langle p_0,s_0 \rangle$.
\For {$k=1,2,\ldots$}
    \State $x_k = x_{k-1} + a_{k-1} p_{k-1}$
    \State $r_k = r_{k-1} - a_{k-1} s_{k-1}$
    \State $w_k = A r_k$
    \State $\nu_{k} = \langle r_k,r_k \rangle$, $~~b_k = \nu_k / \nu_{k-1}$
    \State $\eta_k = \langle r_k,w_k\rangle$, $~~a_k = \nu_k / (\eta_k - ( b_k / a_{k-1} ) \nu_k )$
    \State $p_k = r_k + b_k p_{k-1}$
    \State $s_k = w_k + b_k s_{k-1}$
\EndFor
\EndProcedure
\end{algorithmic}
\end{algorithm}
\end{minipage}
\end{center}
\vspace{.2in}

Notice that the computation of $\nu_k = \langle r_k , r_k \rangle$ can be overlapped
with that of $w_k = A r_k$.  Alternatively, once $w_k = A r_k$ has been formed, the two
inner products $\nu_k = \langle r_k , r_k \rangle$ and $\eta_k = \langle r_k , w_k \rangle$
can be computed simultaneously.
In exact arithmetic, the additional vector $s_k$ is equal to $A p_k$.

The CGCG algorithm can be written in the form (\ref{CGLan}) in much the same way as the
HSCG algorithm.  For the finite precision computation, the relevant formulas are
\begin{eqnarray*}
r_k & = & r_{k-1} - a_{k-1} s_{k-1} + \delta_{r_k} , \\
p_k & = & r_k + b_k p_{k-1} + \delta_{p_k} , \\
s_k & = & A r_k + b_k s_{k-1} + \delta_{s_k} ,
\end{eqnarray*}
where the roundoff terms $\delta_{r_k}$, $\delta_{p_k}$, and $\delta_{s_k}$ 
can be bounded using (\ref{standard}), very similarly to those in HSCG:
\begin{eqnarray*}
\| \delta_{r_k} \| & \leq & \epsilon~ ( \| r_{k-1} \| + 2 \| a_{k-1} s_{k-1} \| )
+ O( \epsilon^2 ) , \\ 
\| \delta_{p_k} \| & \leq & \epsilon~ ( \| r_k \| + 2 \| b_k p_{k-1} \| ) 
+ O( \epsilon^2 ) , \\
\| \delta_{s_k} \| & \leq & \epsilon~ ( \| A r_k \| + 2 \| b_k s_{k-1} \| +
c \| A \|~\| r_k \| ) + O( \epsilon^2 ) .
\end{eqnarray*}
Again, we will assume that the coefficients $a_{k-1}$ and $b_k$ satisfy
the formulas in the CGCG algorithm, with any errors in evaluating
these formulas included in the $\delta_{r_k}$, $\delta_{p_k}$, and
$\delta_{s_k}$ terms.

Eliminating the $s_k$'s and $p_k$'s, we can obtain a three-term recurrence for $r_k$:
\begin{eqnarray*}
r_k & = & r_{k-1} - a_{k-1} ( A r_{k-1} + b_{k-1} s_{k-2} + \delta_{s_{k-1}} ) + \delta_{r_k} \\
    & = & r_{k-1} - a_{k-1} A r_{k-1} - \frac{a_{k-1} b_{k-1}}{a_{k-2}} ( r_{k-2} - r_{k-1} +
\delta_{r_{k-1}} ) - a_{k-1} \delta_{s_{k-1}} + \delta_{r_k} \\
 & = & r_{k-1} - a_{k-1} A r_{k-1} - \frac{a_{k-1} b_{k-1}}{a_{k-2}} ( r_{k-2} - r_{k-1} ) - 
\gamma_k , \\
\end{eqnarray*} 
where
\[
\gamma_k = \frac{a_{k-1} b_{k-1}}{a_{k-2}} \delta_{r_{k-1}} + a_{k-1} \delta_{s_{k-1}} -
\delta_{r_k} .
\]
Defining $q_{k+1} := (-1 )^k r_k / \| r_k \|$, and proceeding exactly as was done for HSCG,
we obtain equation (\ref{HSCGLan1}), where now the roundoff term is
\begin{equation}
\frac{(-1 )^k}{a_{k-1} \| r_{k-1} \|} \gamma_k =
\frac{( -1 )^k}{\| r_{k-1} \|} \left( \frac{b_{k-1}}{a_{k-2}} \delta_{r_{k-1}} + 
\delta_{s_{k-1}} - \frac{1}{a_{k-1}} \delta_{r_k} \right) . \label{fkCGCG}
\end{equation}
Again, this involves only local rounding errors. 
Comparing (\ref{fkHSCG}) and (\ref{fkCGCG}), we see that they look very similar
except that $A \delta_{p_{k-1}}$ in (\ref{fkHSCG}) is replaced by $\delta_{s_{k-1}}$ in (\ref{fkCGCG}), 
where, in exact arithmetic, $s_{k-1} = A p_{k-1}$.

\subsection{GVCG}\label{subsection:GVCG}

This algorithm, developed by Ghysels and Vanroose \cite{GhyVan} and also known as pipelined CG, 
offers the most opportunity for overlapping parallel computations:

\begin{center}
\begin{minipage}{.9\linewidth}
\begin{algorithm}[H]
\caption{Ghysels and Vanroose Conjugate Gradient}\label{gvcg_alg}
\fontsize{10}{10}\selectfont
\begin{algorithmic}
\Procedure{GVCG}{$A$, $b$, $x_0$};~~~{\bf Output}($x_k$)
    \State \textbf{Set} $r_0 = b-Ax_0$, $\nu_0 = \langle r_0,r_0 \rangle$, $p_0 = r_0$, $s_0 = Ap_0$,
    \State \phantom{\textbf{set}} $w_0=s_0$, $u_0 = Aw_0$, $a_0 = \nu_0 / \langle p_0,s_0 \rangle$.
\For {$k=1,2,\ldots$}
    \State $x_k = x_{k-1} + a_{k-1} p_{k-1}$
    \State $r_k = r_{k-1} - a_{k-1} s_{k-1}$
    \State $w_k = w_{k-1} - a_{k-1} u_{k-1}$
    \State $\nu_{k} = \langle r_k,r_k \rangle$, $~~b_k = \nu_k / \nu_{k-1}$
    \State $\eta_k = \langle r_k,w_k\rangle$, $~~a_k = \nu_k / (\eta_k - ( b_k / a_{k-1} ) \nu_k )$
    \State $t_k = Aw_k$
    \State $p_k = r_k + b_k p_{k-1}$
    \State $s_k = w_k + b_k s_{k-1}$
    \State $u_k = t_k + b_k u_{k-1}$
\EndFor
\EndProcedure
\end{algorithmic}
\end{algorithm}
\end{minipage}
\end{center}
\vspace{.2in}

Note that the computation of both inner products $\langle r_k , r_k \rangle$ and
$\langle w_k , r_k \rangle$ required at each iteration can be overlapped with
each other and with the matrix vector product, $t_k = A w_k$, as well as with
some of the vector operations.  In exact arithmetic, the auxiliary vectors 
satisfy $s_k = A p_k$, $w_k = A r_k$, $u_k = A s_k = A^2 p_k$, $t_k = A w_k = A^2 r_k$.

In finite precision arithmetic, the vectors in the GVCG algorithm satisfy
\begin{eqnarray*}
r_k & = & r_{k-1} - a_{k-1} s_{k-1} + \delta_{r_k} , \\
w_k & = & w_{k-1} - a_{k-1} u_{k-1} + \delta_{w_k} , \\
p_k & = & r_k + b_k p_{k-1} + \delta_{p_k} , \\
s_k & = & w_k + b_k s_{k-1} + \delta_{s_k} , \\
u_k & = & A w_k + b_k u_{k-1} + \delta_{u_k} , 
\end{eqnarray*}
where, again using (\ref{standard}), we see that the roundoff terms satisfy
\begin{eqnarray*}
\| \delta_{r_k} \| & \leq & \epsilon~ ( \| r_k \| + 2 \| a_{k-1} s_{k-1} \| )
+ O( \epsilon^2 ) , \\
\| \delta_{w_k} \| & \leq & \epsilon~ ( \| w_{k-1} \| + 2 \| a_{k-1} u_{k-1} \| )
+ O( \epsilon^2 ) , \\
\| \delta_{p_k} \| & \leq & \epsilon~ ( \| r_k \| + 2 \| b_k p_{k-1} \| )
+ O( \epsilon^2 ) , \\
\| \delta_{s_k} \| & \leq & \epsilon~ ( \| w_k \| + 2 \| b_k s_{k-1} \| ) 
+ O( \epsilon^2 ) , \\
\| \delta_{u_k} \| & \leq & \epsilon~ ( \| A w_k \| + 2 \| b_k u_{k-1} \| +
c \| A \|~\| w_k \| ) + O( \epsilon^2 ) .
\end{eqnarray*}

When we try to eliminate auxiliary vectors and form a three-term recurrence for $r_k$,
we find that
\begin{eqnarray*}
r_k & = & r_{k-1} - a_{k-1} ( w_{k-1} + b_{k-1} s_{k-2} + \delta_{s_{k-1}} ) + \delta_{r_k} \\
    & = & r_{k-1} - a_{k-1} w_{k-1} - \frac{a_{k-1} b_{k-1}}{a_{k-2}} ( r_{k-2} - r_{k-1} +
\delta_{r_{k-1}} ) - a_{k-1} \delta_{s_{k-1}} + \delta_{r_k} .
\end{eqnarray*}
It was noted that $w_{k-1} = A r_{k-1}$ in exact arithmetic, so we can write this
recurrence in the form
\begin{eqnarray}
r_k & = & r_{k-1} - a_{k-1} A r_{k-1} - \frac{a_{k-1} b_{k-1}}{a_{k-2}} ( r_{k-2} - r_{k-1} ) -
\nonumber \\ 
 & &  \frac{a_{k-1} b_{k-1}}{a_{k-2}} \delta_{r_{k-1}} - a_{k-1} \delta_{s_{k-1}} + \delta_{r_k} 
 - a_{k-1} ( w_{k-1} - A r_{k-1} ) . \label{rkrecurrence}
\end{eqnarray}
The amount by which the computed vector $r_k$ fails to satisfy a three-term recurrence
now depends not only on local rounding errors, but also on the amount by which
$w_{k-1}$ differs from $A r_{k-1}$.  This will involve rounding errors made at
all previous steps.  To see the size of this difference, subtract $A$ times the 
equation for $r_{k-1}$ from the equation for $w_{k-1}$:
\[
w_{k-1} - A r_{k-1} = w_{k-2} - A r_{k-2} - a_{k-2} ( u_{k-2} - A s_{k-2} ) + \delta_{w_{k-1}}
- A \delta_{r_{k-1}} ,
\]
and apply this recursively to obtain
\begin{equation}
w_{k-1} - A r_{k-1} = w_0 - A r_0 - \sum_{j=0}^{k-2} a_j ( u_j - A s_j ) + \sum_{j=1}^{k-1}
\delta_{w_{j}} - A \sum_{j=1}^{k-1} \delta_{r_j} . \label{wmAr}
\end{equation}
To determine the size of the difference between $u_j$ and $A s_j$, 
subtract $A$ times the equation for $s_{j}$ from that for $u_{j}$ and apply recursively
to find
\begin{eqnarray*}
u_j - A s_j & = & b_j ( u_{j-1} - A s_{j-1} ) + \delta_{u_{j}} - A \delta_{s_j} \\
            & = & b_j b_{j-1} ( u_{j-2} - A s_{j-2} ) + b_j ( \delta_{u_{j-1}} - 
A \delta_{s_{j-1}} ) + ( \delta_{u_j} - A \delta_{s_j} ) \\
            & \vdots & \\
 & = & \left( \prod_{\ell = 1}^{j} b_{\ell} \right) ( u_0 - A s_0 ) + 
\sum_{m=1}^{j} \left( \prod_{\ell = 0}^{m-2} b_{j- \ell} \right) ( \delta_{u_{j-m+1}} -
A \delta_{s_{j-m+1}} ) .
\end{eqnarray*}
Finally, noting that $b_{\ell} = \| r_{\ell} \|^2 / \| r_{\ell -1} \|^2$, one can replace
the above products to obtain
\[
u_j - A s_j = \frac{\| r_j \|^2}{\| r_0 \|^2} ( u_0 - A s_0 ) + \sum_{m=1}^j
\frac{\| r_j \|^2}{\| r_{j-m+1} \|^2} ( \delta_{u_{j-m+1}} - A \delta_{s_{j-m+1}} ) .
\]
Substituting this expression into (\ref{wmAr}) and the result into (\ref{rkrecurrence}),
we see the amount by which $r_k$ may fail to satisfy the three-term recurrence 
that it satisfied in the other algorithms to within {\em local} roundoff errors:  
\begin{eqnarray*}
r_k & = & r_{k-1} - a_{k-1} A r_{k-1} - \frac{a_{k-1} b_{k-1}}{a_{k-2}} ( r_{k-2} - r_{k-1} ) - 
\frac{a_{k-1} b_{k-1}}{a_{k-2}} \delta_{r_{k-1}} - a_{k-1} \delta_{s_{k-1}} + \delta_{r_k}  \\ 
 & & - a_{k-1} \left[ w_0 - A r_0 + \sum_{j=1}^{k-1} \delta_{w_j} - A \sum_{j-1}^{k-1} 
\delta_{r_j} \right. - \\
 & & \left. \sum_{j=0}^{k-2} a_j \left[
\left( \prod_{\ell = 1}^{j} b_{\ell} \right) ( u_0 - A s_0 ) + 
\sum_{m=1}^{j} \left( \prod_{\ell = 0}^{m-2} b_{j- \ell} \right) ( \delta_{u_{j-m+1}} -
A \delta_{s_{j-m+1}} ) \right] \right] .
\end{eqnarray*}
This suggests that the matrix $F_J$ in (\ref{FPLan}) may be significantly larger for
this algorithm than for the others, since it involves roundoff terms
from all steps of the computation, and roundoff terms that are small compared to, say, $\| r_0 \|$
might not be so small compared to $\| r_k \|$.

Note that variants of the rounding error expressions derived in this
subsection have been presented in \cite{Carson,Cools}.  
Note also that a shifted version 
of GVCG \cite{CVnew} has recently been 
proposed to avoid some of the {\em accuracy} problems with the original 
version.  In fact, on all of the test problems in Figure 2 of this paper, 
this version (using an appropriate shift based on estimates of the 
largest and smallest eigenvalues of the matrix) achieves the same ultimate 
level of accuracy as HSCG and CGCG.  However, it is still the case that 
the entries of the matrix $F_J$ in (\ref{FPLan}) are larger than for 
HSCG and CGCG.  Using this version, we have not been able to improve 
significantly on the convergence rate {\em before} the ultimate level of 
accuracy is achieved for the problems shown in Figures 1 and 6.  
Other pipelined CG variants have also been proposed in \cite{CCVnew}.  
The analysis of these procedures is beyond the scope of this paper, but 
again the goal has been improvement in the ultimate level of accuracy.  
In fact, the authors
of \cite{CCVnew} state that ``... it is well known that Krylov subspace
methods may also suffer from delay of convergence due to loss of basis 
orthogonality.... Analyzing the stability issues related to loss of
orthogonality deserves to be treated as part of future work.''

\section{Some Test Problems}

\begin{figure}[ht]
\centerline{\epsfig{file=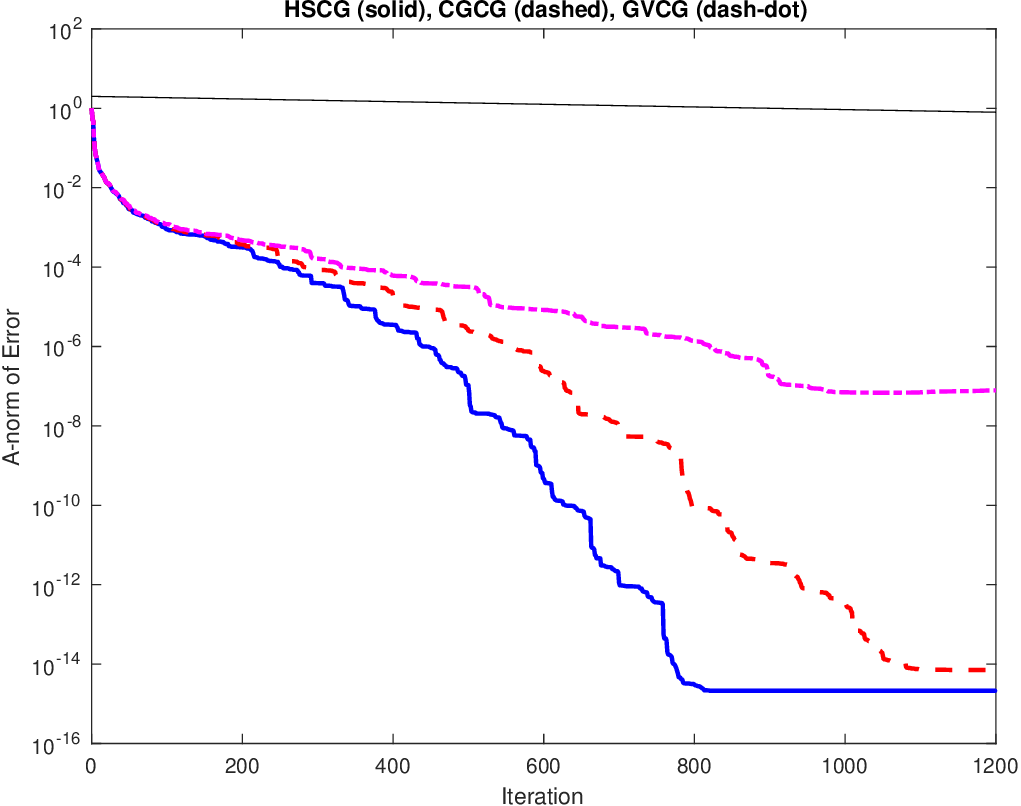,width=3in}}
\caption{Behavior of HSCG (bottom, solid line), CGCG (dashed), and GVCG (dash-dot)
    on the {\tt bcsstk03} matrix.  In exact arithmetic, the exact solution would be
obtained after 112 steps.  The top solid line is the bound (\ref{FPAnorm}).}
\label{fig:bcsstk03}
\end{figure}

The following problem --\verb+bcsstk03+ from the BCSSTRUC1 (BCS Structural Engineering 
Matrices) in the Harwell-Boeing collection \cite{HarwellBoeing} -- was studied in \cite{Carson}.
It is a $112 \times 112$ matrix with condition number $6.8e+6$.  For convenience, we
normalized the matrix so that the matrix we used had spectral norm 1.  In exact arithmetic,
the CG algorithm would obtain the exact solution in at most $112$ steps.  Results of
running HSCG, CGCG, and GVCG are plotted in Figure \ref{fig:bcsstk03}.
We set a random solution vector $x$ and computed $b = Ax$, and we used a zero initial guess $x_0$.
Computations were carried out in MATLAB, using standard double precision arithmetic.
To be sure that we had the true solution $A^{-1} b$ to double precision
accuracy, we solved the linear system $Ax=b$ directly using higher precision and replaced
the original vector $x$ by the one obtained by rounding the multiprecision solution to double precision. 
The figure shows the $A$-norm of the error at each step $k$,
$\langle A^{-1} b - x_k , b - A x_k \rangle^{1/2}$, divided by the $A$-norm of the 
initial error.  Also shown is the upper bound (\ref{FPAnorm}), which is a large
overestimate for all variants.  Note that the different variants of CG not only 
reach different levels of accuracy, but even before the ultimate accuracy level 
is reached, they converge at different rates.
The fastest (in terms of number of iterations) is HSCG, followed by CGCG, with GVCG 
requiring the most iterations.

\begin{figure}[ht]
\centerline{\epsfig{file=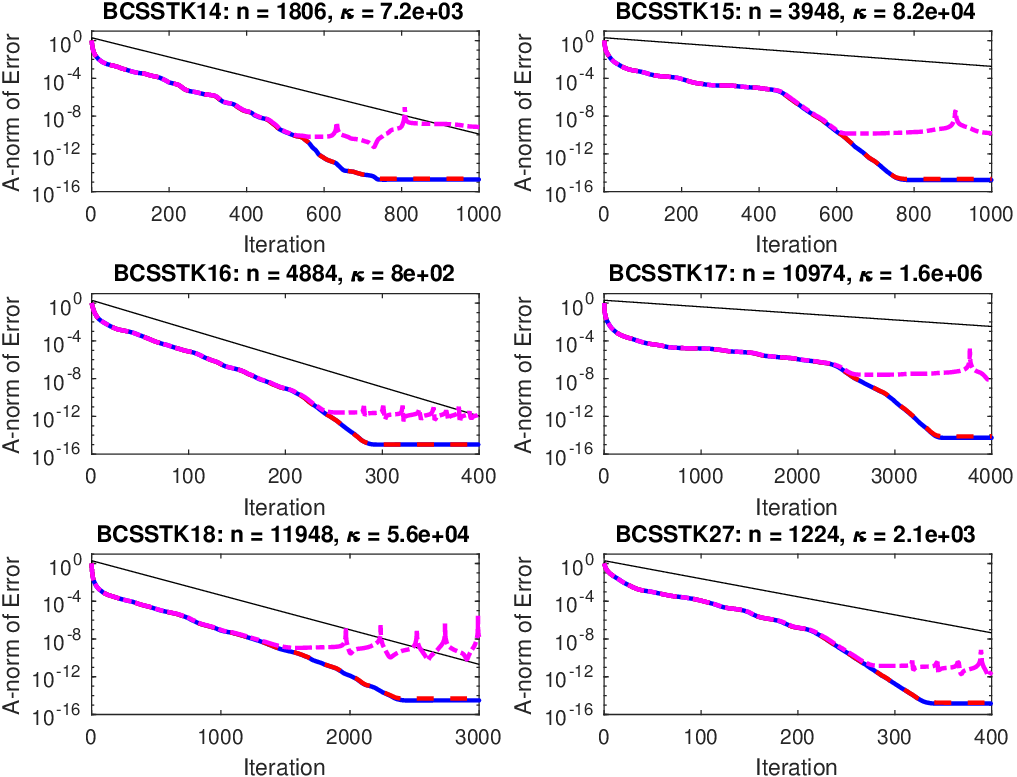,width=3in}}
\caption{Behavior of HSCG (thick solid line), CGCG (dashed, on top of thick solid line), and GVCG (dash-dot)
on matrices from the BCS Structural Engineering Section of the Harwell-Boeing
collection \cite{HarwellBoeing}.  The top solid line is the bound (\ref{FPAnorm}).}
\label{fig:bcsstkall}
\end{figure}

The situation is different, however, for other matrices in this collection.
Figure \ref{fig:bcsstkall} shows the convergence
of HSCG, CGCG, and GVCG for six other test matrices.  Here we used the diagonal
of the matrix as a preconditioner but, to avoid possibly different rounding errors in
preconditioned variants, each matrix was prescaled by its diagonal (that is,
$A$ was replaced by $D^{-1/2} A D^{-1/2}$, where $D = \mbox{diag}(A)$)
and the value $\kappa$ printed on each plot is the condition number of the prescaled matrix.  
Again, we set a random solution vector (and computed the right-hand side
as the product of the prescaled matrix times the random solution vector) and a zero initial guess.  While there
is still some difference in the attainable level of accuracy for the different
variants, until this level is reached, all methods converge at essentially
the same rate. The bound (\ref{FPAnorm}) is shown as well, and while this provides
a good estimate of the actual convergence rate for some of the problems, it is
a large overestimate for others.
Note that here we are plotting the true $A$-norm of the error
$\langle A^{-1} b - x_k , b - A x_k \rangle^{1/2}$ and once $b - A x_k$ starts to
differ substantially from the updated vector $r_k$, one would not expect this quantity to continue
to decrease; we see that with GVCG it may even grow.

For several of these problems, we computed 
\begin{equation}
\max_k \| x_k - ( x_0 + Q_k T_k^{-1} \beta \xi_1 ) \| / \| A^{-1} b \| ,
\label{eps3}
\end{equation}
where $Q_k$ and $T_k$ were computed from the CG residuals and coefficients,
to determine if the tridiagonal linear systems were being solved accurately.  
In all cases, this quantity was tiny, indicating that all variants
are accurately solving the tridiagonal system.  The difference must be in the tridiagonal matrices
that they are producing.  
For all of the problems, the tridiagonal
matrix produced at the last step shown in Figure \ref{fig:bcsstkall} for
the GVCG computation was {\em indefinite}, while the condition numbers of
the final tridiagonal matrices produced in the HSCG and CGCG computations were 
essentially equal to that of $A$.  Note, however, that these computations were run well past
the point where the true and updated residual vectors started to differ significantly, and
the tridiagonal matrices produced while the true and updated residual
vectors were still in close agreement in the GVCG algorithm were 
positive definite and had condition numbers approaching that of $A$, similar to the other
methods.  

We also computed the quantities
\begin{equation}
\epsilon_1 := \max_k \| f_k \| / \| A \| ~ \mbox{ and }~
\epsilon_2 := \max_k | \langle \beta_k q_{k+1} , q_k \rangle | / \| A \| ,
\label{eps12}
\end{equation}
where $f_k$ is the $k$th column of $F_J$ in (\ref{FPLan}), when $T_J$, $Q_J$, 
$\beta_J$, and $q_{J+1}$ come from the CG residuals and coefficients,
to see if conditions (\ref{FJnorm}) and (\ref{successiveinnerprod}) hold.
In all cases $\epsilon_2$ was tiny, that is, a modest multiple of the machine precision.
The same was true for $\epsilon_1$ in the HSCG and CGCG computations, but {\em not} in GVCG.
This might be expected based on arguments in the previous section.  
We observed, however, that for most steps before $r_k$ and $b-A x_k$ started to
differ greatly and before GVCG's ultimate level of accuracy was
reached, the value of $\epsilon_1$ in GVCG, while larger than that in HSCG
and CGCG, was still less than about $1.0e-7$.
We reasoned that for these steps, while the HSCG and CGCG computations behaved like exact
CG for a problem with eigenvalues throughout tiny intervals about the eigenvalues of $A$, 
the GVCG computation might behave like exact CG for a problem with eigenvalues throughout small, 
but not {\em as} small, intervals about the eigenvalues of $A$.  If the interval size made a significant 
difference in the convergence of exact CG, then one would expect slower convergence from GVCG (as seen with \verb+bcsstk03+), 
while if the interval size made little difference in the convergence of exact CG, then one would
expect all three variants to converge at about the same rate (as seen with the other \verb+bcsstk+ problems).

\begin{figure}[ht]
\centerline{\epsfig{file=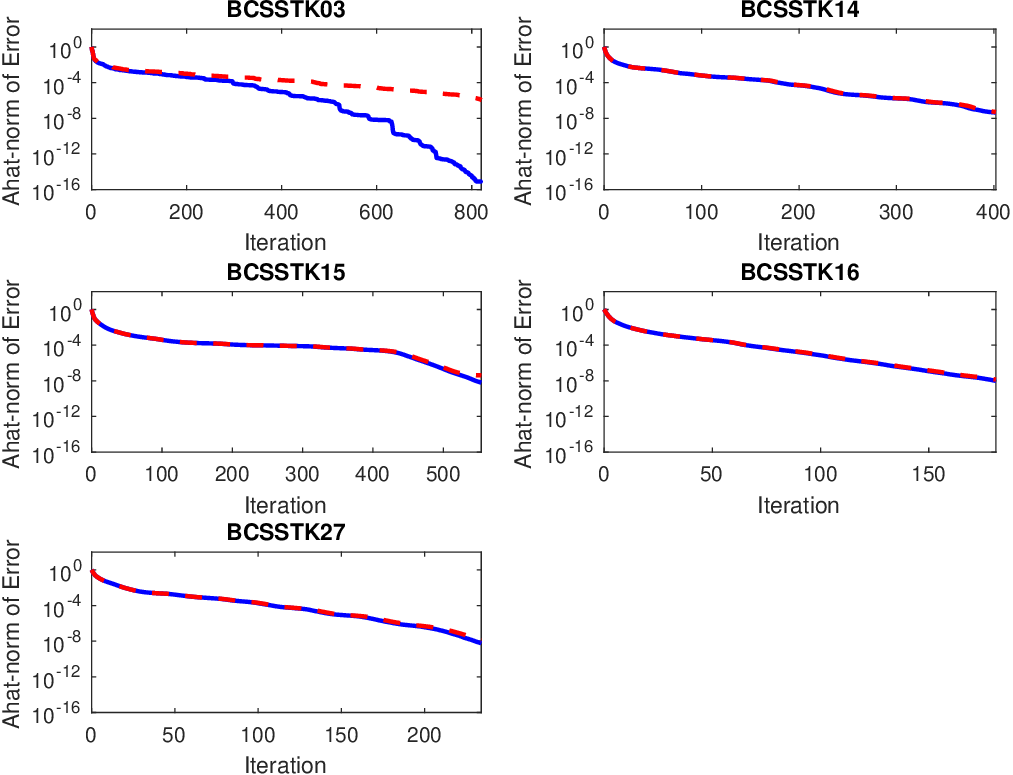,width=4in}}
\caption{Behavior of exact CG for problems with many eigenvalues distributed
throughout intervals of width $1.0e-14$ (solid) and $1.0e-7$ (dashed) about
the eigenvalues of the {\tt bcsstk} matrices.  Results are shown only for
the steps at which $\epsilon_1$ in GVCG remained less than $1.0e-7$.} 
\label{fig:exactCG}
\end{figure}

To test this hypothesis, we computed the eigenvalues of several of the matrices:
\verb+bcsstk03+, \verb+bcsstk14+, \verb+bcsstk15+, \verb+bcsstk16+, and \verb+bcsstk27+.  For each, we formed a larger (diagonal) matrix $\hat{A}$ 
with $11$ eigenvalues distributed about each eigenvalue of the given matrix, in intervals of width
$1.0e-14$ or $1.0e-7$.  Our aim was to determine if the convergence of exact CG
was significantly affected by this interval size, so we used multiple precision
arithmetic and full reorthogonalization of the CG residuals to emulate
exact arithmetic.  We used the same random solution vector for both interval sizes,
and a zero initial guess was used.  The results are shown in Figure \ref{fig:exactCG}.
We ran the computations only to the step where $\epsilon_1$ 
in GVCG started to exceed $1.0e-7$.  As noted above, after this point, the true and
updated GVCG residual vectors started to diverge and progress soon ceased.
The interval sizes $1.0e-14$ and $1.0e-7$ roughly 
coincided with the size of $\epsilon_1$ for HSCG/CGCG and for GVCG, respectively.

Note that the interval width makes a large difference in the convergence of exact CG 
for the matrix $\hat{A}$ associated with the \verb+bcsstk03+ matrix, and the different
CG variants behaved very differently in finite precision arithmetic.  For the other
problems, the interval width makes little difference in the convergence of exact CG
for the matrix $\hat{A}$, and all CG variants behaved similarly in finite precision 
arithmetic (at least, over the steps shown in Figure \ref{fig:exactCG}).
Again, it is our goal to explain the behavior of the different
CG variants at these steps and not at later steps when agreement between true
and updated residuals has been lost.  

\begin{figure}[ht]
\centerline{\epsfig{file=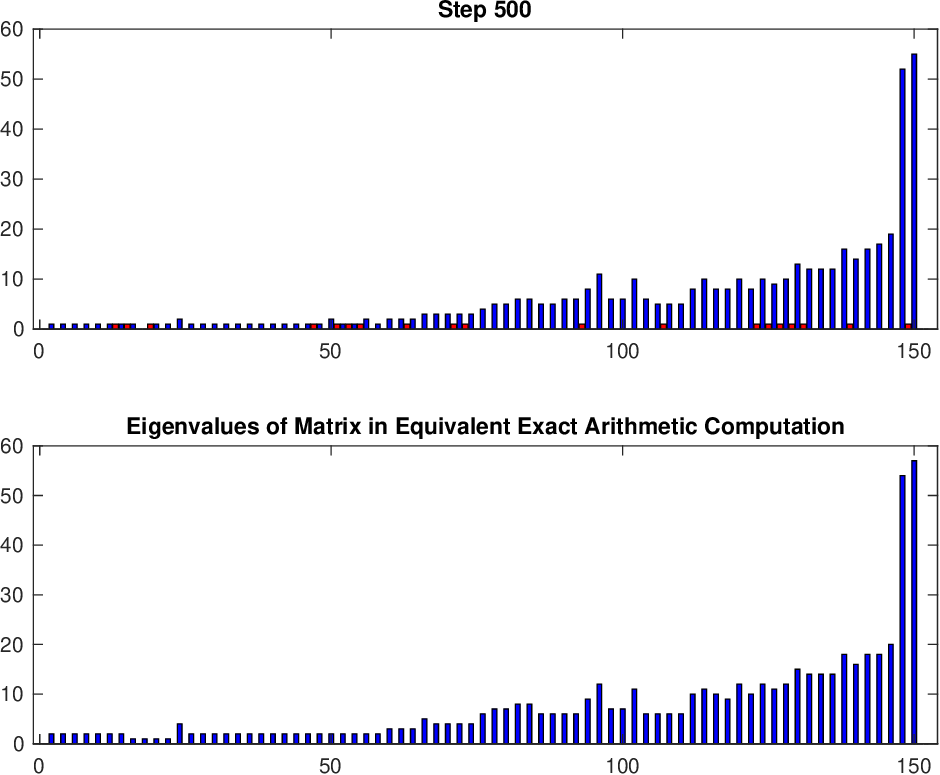,width=4in}}
\caption{Eigenvalues of the tridiagonal matrix produced after 500 steps of HSCG for the 
bcsstk03 problem (top) and eigenvalues of a matrix $T$ for which exact CG
behaves like the finite precision HSCG computation for the first 500 steps (bottom).  Even numbered
bins count eigenvalues within $1.0e-8$ of an eigenvalue of the bcsstk03 matrix and
odd numbered bins count the others. Eigenvalues in odd numbered bins are shown in red
and those in even numbered bins are in blue.
Note that all eigenvalues of $T$ are in even numbered bins.} 
\label{fig:hscgeigs}
\end{figure}

To further understand the behavior of the different algorithms for the \verb+bcsstk03+
problem, we used the procedure described in \cite{Greenbaum89} and
outlined in the Appendix of this paper to construct a matrix
for which exact CG would behave like each of the different finite precision variants for the first
500 steps, which is somewhat past the point where the convergence curves in Figure \ref{fig:bcsstk03}
start to follow different trajectories.  [See the Appendix for a more precise
description of the relation between the finite precision computations and
exact CG for this larger matrix.] There are many different matrices $T$
for which exact CG applied to $T$ would match the behavior of each of these 
finite precision computations \cite{Meurant}, and the extension procedure given in
\cite{Greenbaum89} does not necessarily produce the matrix $T$ with eigenvalues
in the smallest possible intervals about the eigenvalues of $A$.  
Still, we found that for HSCG, it produced a matrix with eigenvalues in intervals of
width  $2.2e-9$ about the eigenvalues of $A$, for CGCG the interval width was $1.7e-8$,
and for GVCG it was $3.0e-7$.  The theory is not precise
enough to predict the difference in interval width for HSCG and CGCG based
on the size of $\epsilon_1$, $\epsilon_2$, and $\epsilon_3$, but it does
appear that this interval width must be larger for CGCG, accounting for the
somewhat slower rate of convergence.
For illustration, the eigenvalues of the tridiagonal matrix
produced at step 500 by the HSCG computation and the eigenvalues of the matrix $T$ for which
exact CG matches the behavior of HSCG for the first 500 steps are plotted in Figure \ref{fig:hscgeigs}.
The figure uses histograms with bins of width $1.0e-8$ about each eigenvalue of the \verb+bcsstk03+
matrix and bins to represent values in between.  Thus, if $\lambda_1 < \lambda_2 < \ldots$ are
the distinct eigenvalues of the \verb+bcsstk03+ matrix, then any eigenvalues 
less than $\lambda_1 - 1.0e-8$ are counted in bin 1, any between $\lambda_1 - 1.0e-8$ and
$\lambda_1 + 1.0e-8$ are counted in bin 2, those between $\lambda_1 + 1.0e-8$ and 
$\lambda_2 - 1.0e-8$ are counted in bin 3, etc.  Note that all eigenvalues of the matrix $T$
are in even numbered bins, meaning that they are within $1.0e-8$ of an eigenvalue of the \verb+bcsstk03+ 
matrix, but as noted above, they were actually within $2.2e-9$ of these eigenvalues.

\begin{figure}[ht]
\centerline{\epsfig{file=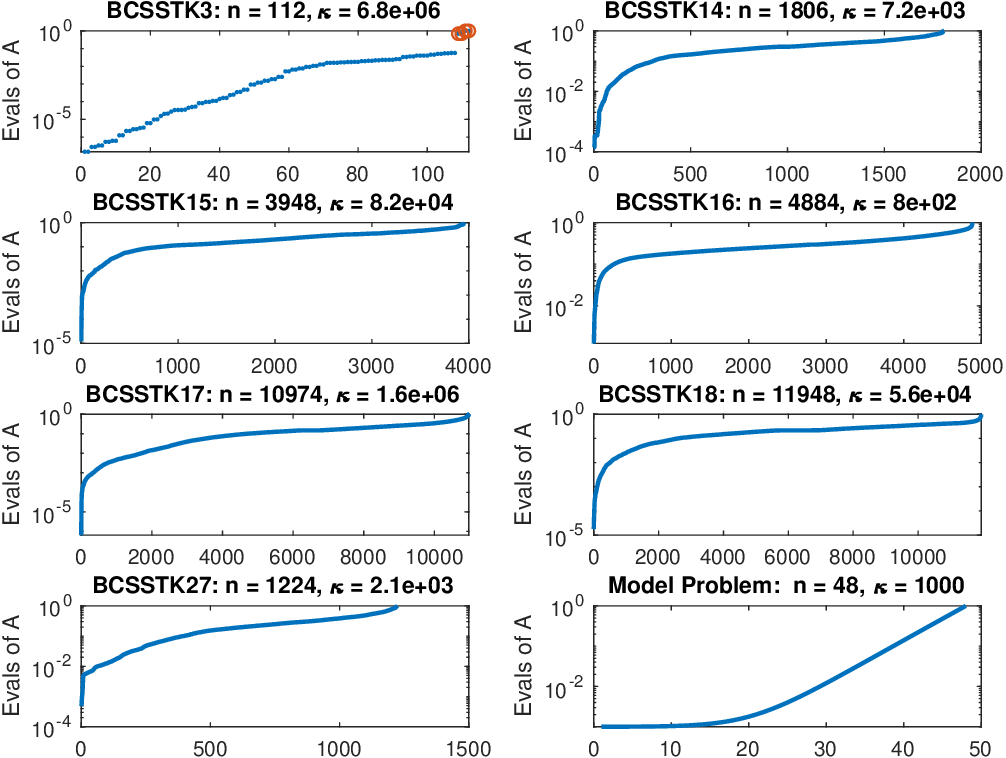,width=4in}}
\caption{Eigenvalues of the bcsstk matrices and of a model problem from \cite{GreenStrak}.} 
\label{fig:evals}
\end{figure}

Finally, to understand what sorts of eigenvalue distributions lead to different
convergence curves for the different finite precision implementations, we plot the
eigenvalues of these matrices on a log scale in Figure \ref{fig:evals}.  
Note that the \verb+bcsstk03+ matrix has four large eigenvalues (marked with $\circ$'s in the figure)
that are well-separated from the others.  It is known that exact CG applied to a matrix 
with eigenvalues throughout
intervals around large well-separated eigenvalues will converge significantly slower
than exact CG for a matrix that has just the few discrete outliers, since the CG polynomial
associated with the intervals will have multiple roots within each interval.  Moreover, the
size of the intervals will determine the frequency with which roots are put down in the intervals.
See \cite{Greenbaum89}.  The other \verb+bcsstk+ matrices have more small well-separated eigenvalues, 
which present less of a problem for finite precision computations \cite[Ch. 4]{Greenbaumbook}.  
Moreover, ignoring the very small eigenvalues, the middle eigenvalues of the other \verb+bcsstk+
matrices range over fewer orders of magnitude than the \verb+bcsstk03+ eigenvalues.  

\begin{figure}[ht]
\centerline{\epsfig{file=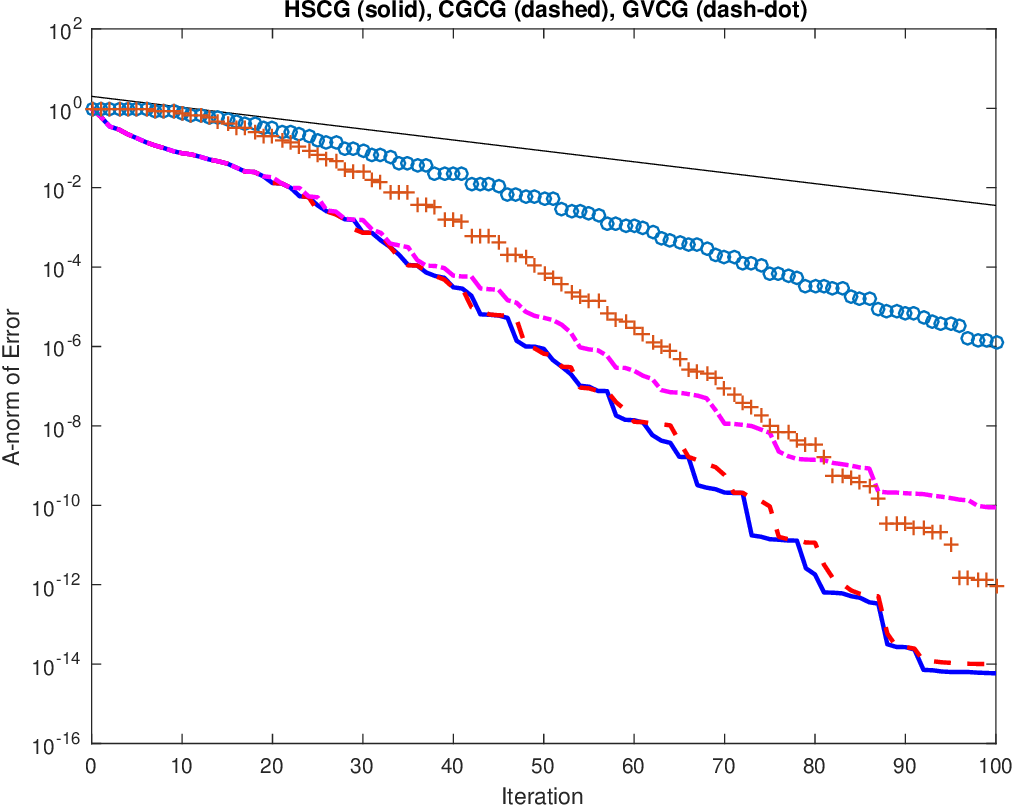,width=3in}}
\caption{Behavior of HSCG (bottom, solid line), CGCG (dashed), and GVCG (dash-dot)
on a model problem with eigenvalues given by (\ref{modeleigs}) for $n = 48$,
$\rho = 0.8$.  In exact arithmetic, the exact solution would be obtained after
$48$ steps.  The top solid line is the bound (\ref{FPAnorm}), the $\circ$'s
are the quantity (\ref{sharpbound}) with $\delta = 1.0e-7$,
and the $+$'s are the quantity (\ref{sharpbound}) with $\delta = 1.0e-14$.} \label{fig:test1}
\end{figure}

For further verification, we include a test problem from \cite{GreenStrak} whose eigenvalues
are highly spread out at the upper end.  Taking $n=48$ and $\rho = 0.8$, we formed a matrix
with the following eigenvalues:
\begin{equation}
\lambda_1 = 0.001,~~\lambda_n = 1 ,~~
\lambda_i = \lambda_1 + \frac{i-1}{n-1} ( \lambda_n - \lambda_1 ) \rho^{n-i} ,~~
i=2, \ldots , n-1 . \label{modeleigs}
\end{equation}
These eigenvalues are also shown in Figure \ref{fig:evals}, in the bottom right subplot.
We chose random orthonormal eigenvectors, a random solution vector, and a zero initial guess for the solution.
Results of running the HSCG, CGCG, and GVCG algorithms are plotted in Figure \ref{fig:test1}.
Also plotted is the upper bound (\ref{FPAnorm}) using $\kappa = 1000$ and the
quantity (\ref{sharpbound}) using $\delta = 1.0e-14$ and using $\delta = 1.0e-7$.
The minimax polynomial on the union of intervals was computed using the Remez algorithm \cite[pp.~280-289]{Blum}.
Note that for this problem, the interval size makes a significant difference in the
size of the minimax polynomial on the union of intervals.

Although the exact solution would be obtained after $48$ steps using exact arithmetic,
all of the finite precision implementations required about $100$ iterations to achieve
their best level of accuracy.  Note, however, that while the convergence curves for HSCG
and CGCG are very similar before this point, that of GVCG is significantly
worse.  Again we looked at the quantity in (\ref{eps3}) and at
$\epsilon_1$ and $\epsilon_2$ in (\ref{eps12}) and found that all were moderate
multiples of the machine precision {\em except} $\epsilon_1$ in GVCG, which was
$7.5e-5$.  Note that the finite precision GVCG computation {\em cannot} be
equivalent to exact CG for a matrix with eigenvalues in intervals of width
$1.0e-14$ about the eigenvalues of $A$ since the GVCG convergence curve goes
above the upper bound (\ref{sharpbound}), which holds for exact CG for all such
matrices.  The GVCG convergence curve does, however, remain below the upper
bound for exact CG applied to matrices with eigenvalues in intervals of width
$1.0e-7$ about the eigenvalues of $A$.

\section{Conclusions}
In order to relate the behavior of a finite precision CG computation to that of
exact CG for a matrix with eigenvalues in small intervals about the eigenvalues of $A$ using
the analysis in \cite{Greenbaum89}, the finite precision computation must satisfy
(\ref{FJnorm}) and (\ref{successiveinnerprod}) for some small number $\epsilon \ll 1$.
Even if $\epsilon$ is not so small, the finite precision computation may satisfy the 
bounds (\ref{FPresnorm}) and (\ref{FPAnorm}) provided only that the tridiagonal matrix
that it produces has its eigenvalues essentially between the largest and smallest
eigenvalues of $A$ and provided formula (\ref{xJLan}) is satisfied to a close approximation.
These are the conditions that one should check in designing new implementations.

It was observed numerically that all three of the implementations considered here --
HSCG, CGCG, and GVCG -- satisfied (\ref{successiveinnerprod}) and (\ref{xJLan}).
The only difference was that GVCG did not satisfy (\ref{FJnorm}) as well as the 
others.  Yet even GVCG satisfied (\ref{FJnorm}) reasonably well up to a point,
and after that point progress soon ceased.  What happened before that point?  For most
problems, failure to satisfy (\ref{FJnorm}) to the order of machine precision made 
little difference -- the GVCG computation behaved like exact CG for a matrix with
eigenvalues in small but not tiny intervals about the eigenvalues of $A$, and the behavior
of exact CG was not very sensitive to the precise size of the intervals.  For the \verb+bcsstk03+
matrix, however, like the model problem with eigenvalues highly spread out at the upper end,
that is not the case.  The behavior of exact CG for a matrix with eigenvalues 
throughout small intervals about the eigenvalues of those matrices is very sensitive to the size
of the intervals, and the different CG variants behaved differently in finite precision
arithmetic even before the final level of accuracy was reached.
Thus, in evaluating a new implementation of the CG algorithm, one
must be careful to include test problems that stress the effect of rounding errors
on the convergence rate.  It may be that the cost of extra iterations is outweighed
by the advantages of parallelism or other things, but this should be included
in the overall evaluation of the method.

\vspace{.1in}
\noindent
{\bf Acknowledgment:} The authors thank the referees for their careful reading of
this manuscript and their many helpful suggestions.

\section*{Appendix:  Software}

The MATLAB codes used to produce plots in this paper can be found at
\verb+https://github.com/HexuanLiu/Conjugate_gradient+.

The most interesting of these is \verb+extendT.m+, which takes as input a
symmetric positive definite matrix $A$, a symmetric tridiagonal matrix $T_J$
and a set of unit vectors $q_1 , \ldots , q_J$ stored as columns of a matrix $Q_J$
(such as those returned by \verb+HSCG.m+, \verb+CGCG.m+, or \verb+GVCG.m+),
and the number of digits \verb+ndigits+ to use with MATLAB's variable precision
arithmetic toolbox.  It returns a multiprecision symmetric tridiagonal matrix \verb+T_vpa+
that is an extension of $T_J$ whose eigenvalues are, hopefully, all close to
eigenvalues of $A$.  It also returns a multiprecision array \verb+Q_vpa+ whose
columns each have norm 1 and such that \verb+ A * Q_vpa+ is approximately equal to
\verb+ Q_vpa * T_vpa+.  It uses the procedure outlined in \cite{Greenbaum89} to
construct \verb+T_vpa+ and \verb+Q_vpa+.  This procedure is described below.

Assume that the input variables satisfy
\begin{equation}
A Q_J = Q_J T_J + \beta_J q_{J+1} \xi_J^T + F_J . \label{A1}
\end{equation}
Let $T_J S_J = S_J \Theta_J$ be an eigendecomposition of $T_J$, and define
$Y_J := Q_J S_J$.  Multiplying (\ref{A1}) on the right by $S_J$, we have
\begin{equation}
A Y_J = Y_J \Theta_J + \beta_J q_{J+1} \xi_J^T S_J + F_J S_J . \label{A2}
\end{equation}
Let $y_1 , \ldots y_J$ denote the columns of $Y_J$ (referred to as Ritz vectors) and let
$\theta_1 , \ldots , \theta_J$ denote the diagonal entries of $\Theta_J$ (Ritz values).
Define a Ritz value $\theta_i$ to be {\em well-separated} from the others if
\[
\min_{k \neq i} | \theta_k - \theta_i | > \verb+(cluster_width)+ \| A \| ,
\]
where, initially, \verb+cluster_width+ is taken to be the square root of the machine precision;
otherwise, consider it to be part of a cluster.  For clustered Ritz values, define
a cluster vector by
\[
y^C := \frac{1}{w_C} \sum_{\ell \in C} S_{J, \ell} y_{\ell} ,~~
w_C = \left( \sum_{\ell \in C} ( S_{J, \ell} )^2 \right)^{1/2} ,
\]
where $S_{J, \ell}$ is the $(J, \ell )$-entry of $S_J$ and the sum is over all Ritz
vectors corresponding to Ritz values in the cluster.  Define a {\em cluster value} by
\[
\theta_C = \frac{1}{2} \left( \min_{\ell \in C} \theta_{\ell} + \max_{\ell \in C} \theta_{\ell} \right) .
\]

We will say that a Ritz vector $y_i$ corresponding to a well-separated Ritz value is {\em converged}
if $\beta_J | S_{J, i} | \leq \verb+(conv_tol)+ \| A \|$, where, initially, \verb+conv_tol+ is taken
to be the square root of the machine precision; otherwise, it is {\em unconverged}. 
We will say that a cluster vector $y^C$ is {\em converged} if $\beta_J w_C \leq \verb+(conv_tol)+ \| A \|$
and {\em unconverged} otherwise.
Let $\hat{Y}_{m}$ have $m$ columns consisting of the unconverged Ritz vectors and the unconverged cluster vectors.

Assuming that $\epsilon_1$ and $\epsilon_2$ in (\ref{eps12}) are on the order of machine precision, 
it is argued in \cite{Greenbaum89}
that the columns of $\hat{Y}_{m}$ should be almost orthonormal, $q_{J+1}$ should be almost orthogonal
to the columns of $\hat{Y}_{m}$, and $q_J$ should be almost equal to a linear combination of these columns.
Since not all CG variants maintain $\epsilon_1$ and $\epsilon_2$ at the level of machine precision, information
is printed out to show how closely these properties are satisfied, and the user is given
an opportunity to adjust \verb+cluster_width+ and \verb+conv_tol+ to better satisfy these properties.
For HSCG applied to the bcsstk03 problem, we achieved the best results by
taking \verb+cluster_width+ and \verb+conv_tol+ to be $1.0e-9$, while we left them at the square root of
the machine precision for CGCG, and we set them to $1.0e-6$ for GVCG. 

Once the columns of $\hat{Y}_{m}$ are determined, the rest of the code is straightforward.
The next Lanczos vector $q_{J+1}$ is modified (slightly) to be {\em exactly} (that is, to \verb+ndigits+ precision)
orthogonal to the columns of $\hat{Y}_{m}$.  Successive vectors satisfy the usual 3-term Lanczos recurrence, with the
coefficients being used to extend $T_J$, but the recurrence is perturbed to make the new vectors
exactly orthogonal to each other and to the columns of $\hat{Y}_{m}$.  This means that the algorithm will
terminate with $q_{J+n-m}$ equal to $0$ and with a matrix \verb+T_vpa+ of size $J+n-m$ whose eigenvalues
should all be close to eigenvalues of $A$ (assuming, again, that $\epsilon_1$ and
$\epsilon_2$ are sufficiently small).

The driver code, \verb+bcsstk03magic.m+, runs either HSCG, CGCG, or GVCG and then calls \verb+extendT+ to
extend the tridiagonal matrix to one whose eigenvalues are close to eigenvalues of $A$.  It then
runs \verb+cg_vpa+, a variable precision CG code using full reorthogonalization, with the matrix \verb+T_vpa+
and right-hand side equal to the first unit vector. It plots the 2-norm
of the residual and the \verb+T_vpa+-norm of the error at each step
of the multiprecision CG computation to demonstrate that these values
fall right on top of those from the finite precision computation with $A$.
More precisely, the 2-norms of the residuals in exact CG for \verb+T_vpa+
exactly match the 2-norms of the vectors $r_k$ in the finite precision computation,
assuming that the tridiagonal matrix from the finite precision computation
(of which \verb+T_vpa+ is an extension) consists of entries
satisfying the exact coefficient formulas in the algorithms; the match between
the \verb+T_vpa+-norm of the error in exact CG applied to \verb+T_vpa+ and
the quantities $\langle r_k , A^{-1} r_k \rangle$ in the finite precision
computations is very close but not exact.

\end{document}